\documentclass{article}
\usepackage{amssymb, amsmath,amsthm}
\def\R{{\bf R}}

\newtheorem{corollary}{Corollary}
\renewcommand{\epsilon}{\varepsilon}
\renewcommand{\phi}{\varphi}
\newtheorem{remark}{Remark}

\begin{document}
\parindent0 in
\parskip 1 em
\title{On the Discrete Spectrum of Generalized Quantum Tubes}
\author{Christopher Lin and Zhiqin Lu\footnote{
The second
author is partially supported by  NSF Career award DMS-0347033 and the
Alfred P. Sloan Research Fellowship.}
}

\date{February 9th, 2005}

\maketitle
\section{Introduction}\label{section1}
The existence of discrete spectrum of the Laplacian on manifolds is an 
interesting phenomenon in both geometry and physics.  In geometry, discrete 
spectrum characterizes compact manifolds as it is the single constituent for 
their spectrum, while it is often unknown that whether or not a complete 
noncompact manifold has discrete spectrum.  On the other hand, in physics the 
discrete spectrum of the Laplacian represents quantization of energy of a 
free nonrelativistic particle and correspondingly the localization of the 
wave functions.  The physics terminology of bound states are thus the 
eigenfunctions corresponding to points in the discrete spectrum.    

The mathematical problem we study in this paper is the existence of the 
discrete spectrum of the Dirichlet Laplacian on a noncomplete noncompact 
manifold called the quantum tube, which is built as a type of 
tubular neighborhood about an immersed manifold in Euclidean space.  The 
Dirichlet Laplacian is defined as the unique self-adjoint operator 
associated to the closed, positive symmetric quadratic form 

\begin{equation}
Q(\psi, \phi) = \int\langle \nabla\psi, \nabla\phi \rangle \quad 
\forall\psi, \phi \in C_)^\infty, 
\end{equation}    
where we took the metric inner product of the gradients over the 
Riemannian manifold. Let $\Omega$ be a complete manifold with boundary. 
The Dirichlet Laplacian may be thought of as the Laplacian with Dirichlet 
boundary condition. In fact the space 
$\left\{f \, | \, f\in C^{\infty}(\Omega)\,
\text{and} \, f|_{\partial\Omega} = 0 \right\}$ 
is an operator core for the Dirichlet Laplacian.  
Let $\sigma(\Delta)$ and $\sigma_{ess}(\Delta)$ denote the total and 
essential spectrum of the Dirichlet Laplacian respectively.  We 
have the following standard expressions for the bottom of the respective 
spectra on a noncompact manifold $\Omega$:

\begin{equation}\label{infess}
\inf \sigma_{ess}(\Delta) = 
\sup_{K} \inf_{f\in C_{0}^{\infty}( \Omega \setminus K)} 
\frac{\int_\Omega |\nabla f|^2}{\int_\Omega f^2},
\end{equation}
and 
\begin{equation}
\inf\sigma(\Delta) = 
\inf_{f \in C_{0}^{\infty}(\Omega)}
\frac{\int_\Omega |\nabla f|^2}{\int_\Omega f^2},
\end{equation}
where $K$ are compact subsets of $\Omega$.   

As we mentioned, the discrete spectrum of the Laplacian may not exist for 
complete noncompact manifolds.  However, in our previous paper 
\cite{ll-1} we showed that if we were to thicken a complete noncompact 
manifold by one more dimension, then under suitable curvature assumptions 
the discrete spectrum of the Dirichlet Laplacian on the resulting layer 
will be nonempty.  In 
particular, we 
proved that $\inf\sigma(\Delta) < \inf\sigma_{ess}(\Delta)$ for the 
layer.  In \cite{ll-1} the layer we studied is called the quantum layer, due
to the fact that the additional dimension induces quantization of the 
spectrum.  To be more precise, the quantum layer was defined as a 
tubular neighborhood about an orientable, complete noncompact manifold 
$\varSigma$ of dimension $n\geq 2$, 
isometrically immersed in $\mathbb{R}^{n+1}$, with its width $a$ 
suitably controlled by the second fundamental form of the immersion so that 
the tubular neighborhood is also an isometric 
immersion of $\varSigma \times (-a,a)$ into $\mathbb{R}^{n+1}$.  In addition
we required that $\varSigma$ is parabolic (see section~\ref{section2} 
for precise definition) and asymptotically flat.  Then under suitable 
integral-curvature assumptions on $\varSigma$ we can deduce that 
$\inf\sigma(\Delta) < \inf\sigma_{ess}(\Delta)$ for the quantum layer.\footnote{      
The terminology of quantum layer was actually first introduced in 
the paper \cite{dek} by Duclos, Exner, and   
Krej\v{c}i\v{r}\'{\i}k.  The inspiration comes from mesoscopic physics .}  

 In this paper, we will take the generalization further by 
allowing the codimension
of the immersion of $\varSigma$ to be arbitrary.  Whereas in our last paper 
the geometric challenge in enlarging the dimension of the immersed manifold 
$\varSigma$ was overcome by the use of parabolicity, 
the main geometric challenge 
presented here by the arbitrary codimension is to untangle the 
more complicated ways in which the geometry of $\varSigma$ interacts 
with that of the ambient Euclidean space.  In the quantum layer case, 
where the codimension is $1$, the orientability of $\varSigma$ allowed us to 
define an immersion  
$\mathcal{L}: \varSigma\times (-a,a) \longrightarrow \mathbb{R}^{n+1}$ 
by $\mathcal{L}(x,u) = x + uN$, where $N$ is a global unit vector field on 
$\varSigma$.  The existence of $N$ on $\varSigma$ creates a global vertical 
(or fiber) coordinate for the quantum layer.  
For codimension greater than $1$, we no longer have such global fiber 
coordinates.  Therefore to create a $k$-dimensional tubular neighborhood 
about $\varSigma$, we opt for a more intrinsic definition by lifting to 
the normal bundle (see definitions below).  We will use the term quantum 
tube for this more generalized notion of a tubular neighborhood about 
$\varSigma$.  This new terminology brings to mind the image of a tube about 
a curve in space, although this paper does not deal with immersed 
one dimensional 
manifolds (for more remarks on this issue, see next section).  On the quantum 
tube we will perform similar analysis as in 
$\cite{ll-1}$ and prove again that the discrete spectrum is nonempty 
by showing $\inf\sigma(\Delta) < \inf\sigma_{ess}(\Delta)$, in 
particular the ground state (first eigenvalue) 
$\lambda_0 = \inf\sigma(\Delta)$ exists.  However, 
we will see that the quantum tube has new geometric features that make 
this paper worthwile.

First, there is the problem of the nontriviality of the normal bundle of 
$\varSigma$ if the codimension is greater than $1$.  For the quantum 
layer, where the normal bundle of $\varSigma$ was trivial, to study 
the Laplacian on $\varSigma$ we merely needed to concern ourselves 
with the tangent bundle of $\varSigma$ and its connection.  However, 
if the normal bundle is non-trivial, we may need to incorporate the 
normal bundle and the normal connection on $\varSigma$.  In fact, the 
geometry of the quantum tube as revealed by its metric tensor, illustrates 
this point (see section~\ref{section2}).  In the next section we will see that 
under our choice of a  coordinate system, 
the metric tensor of the quantum tube will contain normal 
covariant derivatives of the fiber coordinate vectors along the horizontal 
coordinate frames.  The metric tensor will be in 
non-block form as a result, which presents a difficulty in evaluating and 
estimating the volume measure on the quantum tube.  
One way to bypass this difficulty is to assume that these normal covariant 
derivatives vanish, which corresponds to the existence of parallel frames 
on $\varSigma$ with respect to the normal connection.  However, 
we know that on any vector bundle parallel frames exist 
if and only if the connection is flat, i.e., with 
identically zero curvature (see equation (1.1.4), Proposition 1.1.5, and 
Theorem 1.4.3 in \cite{pt}).  To 
Assume that the immersion of $\varSigma$ has flat normal bundle is a 
restrictive condition, even if the ambient manifold is Euclidean space.  
However, as we will see in sections~\ref{section2} and~\ref{section5}, 
we do not need the assumption of flat normal bundle and the ensuing 
technical difficulty mentioned above could be resolved with some work.  

The second  feature of this paper is the further generalization of 
the curvature assumptions in \cite{dek} and \cite{ll-1}.  Due 
to the arbitrary dimension  and codimension of $\varSigma$, we had to consider 
the average of  even elementary symmetric functions $K_{2p}$
of the second fundamental form.  
These isometric invariants $K_{2p}$ turn out to be quite well-known in the literature.  
We were informed  by Mazzeo that these are 
the famous tube invariants of Weyl. A 
comprehensive book by 
Gray\cite{gray} discusses these tube invariants and their generalizations 
in much detail.  

Next, we give a brief review of basic notions in isometric immersion.  
Let $TN$ denote the tangent bundle of a manifold $N$, and $T_{p}N$ its 
fiber at $p\in N$.  
Given an isometric immersion $\Sigma \hookrightarrow N$, 
 we shall identify the image of the immersion
 with the manifold $\Sigma$ itself. The Riemannian metric 
on $N$ determines a sub-vector bundle of $TN$ called the normal bundle 
of $\Sigma$ which we denote by $T^{\bot}\Sigma$.  At each point its fiber satisfies
$T_{p}^{\bot}\Sigma \oplus T_{p}\Sigma = T_{p}N$.  Formally we may 
define the normal bundle as the quotient bundle 
$T^{\bot}\Sigma = TN\diagup T\Sigma$. 
For the sake
of convenience, we will denote tangent
vector fields on $\Sigma$ by Roman letters and normal vector fields on $\Sigma$
by Greek letters.
 Let $\nabla$ be the Riemannnian
connection on the ambient manifold $N$.  Then the
normal connection $\nabla^{\bot}$ on the normal bundle is 
defined by
$\nabla^{\bot}_X \eta = (\nabla_X \eta)^{\bot}$.  
From the normal connection we also get the normal curvature 
tensor $R^{\bot}$ defined by 
$R^{\bot}(X,Y)\eta = \nabla^{\bot}_Y\nabla^{\bot}_X\eta 
- \nabla^{\bot}_X\nabla^{\bot}_Y\eta 
+ \nabla^{\bot}_{[X,Y]}\eta$. 

We denote the second fundamental
form of the immersion by $\vec{A}$, 
which is defined by the symmetric bilinear tensor
$\vec{A}(X,Y) = (\nabla_{X}Y)^{\bot}$.  
The shape operator associated with a 
normal vector field $\eta$ on $N$ is the self-adjoint operator  
$S_{\eta}(X) = -(\nabla_{X}\eta)^{T}$, which is related to 
the second fundamental form by 
$\langle \vec{A}(X,Y), \eta\rangle = 
\langle S_{\eta}(X), Y\rangle$.  
Moreover, since the shape operator is a tensor we can 
define the elementary symmetric functions $C_j(S_\eta)$, 
$j = 1, 2,\cdots, n$, as genuine smooth functions on $\Sigma$ and 
compute them locally.          

Now, suppose the ambient manifold is the Euclidean space 
$\mathbb{R}^{m}$.  
  Then it makes sense to define the map 
$f: T^{\bot}\Sigma \longrightarrow \mathbb{R}^{m}$ given by  
\begin{equation}\label{E:imm}
f(x,\xi) = x + \xi.
\end{equation}
The map~\eqref{E:imm} is well defined because every 
$\xi \in T_{x}^{\bot}\Sigma$ is 
naturally identified with a unique vector in $\mathbb{R}^{m}$ by 
  $\xi = a_i\frac{\partial}{\partial x_i}$,
 where $\{x_1,\cdots, x_{n+k}\}$ are the standard coordinates of 
$\mathbb{R}^{n+k}$.  
$f$ is also clearly smooth.  
Moreover, with the Euclidean metric $ds_{E}^{2}$ we can compute the
length $\|\xi\|$ of any normal 
vector $\xi$ on $\Sigma$, and we can consider 
the bilinear norm $\|\vec{A}\|$ - 
which becomes a smooth function on $\Sigma$.    

Let us introduce the following conditions on a complete manifold $\Sigma$:

A1) $\|\vec{A}\| \leq \epsilon_{o}$ 
         for some $\epsilon_{o} > 0$

A2) $\|\vec{A}\|(x) \rightarrow 0$ 
         as $d(x, x_0) \rightarrow \infty$ for
         any fixed point $x_0 \in \Sigma$.

\vskip 0.6cm 

\newtheorem{definition}{Definition}
\begin{definition}
Let $\varSigma \hookrightarrow \mathbb{R}^{n+k}$ be an 
 isometric immersed,
complete, noncompact  n-dimensional oriented manifold 
with $n \geq 2$  satisfying condition A1).  Let $F$ be the submanifold 
of $T^{\bot}\varSigma$ defined by 
$F = \{(x,\xi)\in T^{\bot}\varSigma \hskip 0.1cm | 
\hskip 0.1cm \|\xi\| < r\}$,     
where r is any fixed positive number satifying 
$r \leq (\sqrt{k}\hskip 0.05cm \epsilon_{o})^{-1}$. 
Then we define an order-$k$ quantum tube with radius r as the 
Riemannian manifold $(F, f^*(ds_{E}^2))$, 
where $f^*(ds_{E}^{2})$ is the pullback metric induced by 
the map $f$ given in (1) restricted to F.  We
call $\varSigma$ the base manifold of the quantum tube $F$. 
\end{definition}

Our definition above requires that $f$ restricted to $F$ be an 
immersion 
into $\mathbb{R}^{n+k}$.  The condition that guarantees this is exactly
given by the requirement on the radius $r$, as we will see 
in section~\ref{section2}.  Note that we do not require that 
$F$ be embedded in Euclidean space, thus we allow self-intersections
in $\mathbb{R}^{n+k}$, and the spectral results we derive still hold.  
It is also clear that our definition 
generalizes the definition given in 
\cite{ll-1} of the quantum layer.  The definition we give 
above will lend itself to features in 
the geometry of the quantum tube $F$ that was not observed previously 
in the quantum layer of \cite{ll-1}, as remarked earlier. 

Having formally defined the quantum tube, we would like to pause to remark 
on how its peculiar structure allows us to study its spectrum.  The quantum 
tube is clearly a noncompact, noncomplete manifold, and contrary to 
compact manifolds (with or without boundary) and complete manifolds, 
there are no standard techniques to study the Laplacian on such manifolds.  
However, the quantum tube is constructed so that there is the base manifold 
(the zero-section in the normal bundle) 
as the horizontal part and the fibers as the vertical part.  In particular 
this allows us to construct test functions in the form of a product of a 
function on the base manifold and a function on the fibers.  Therefore, 
we can more or less use special properties of the horizontal 
and vertical parts of the test function separately in the same analysis.  
  As we will see 
in section \ref{section5}, this allows us 
use special functions associated to the parabolicity on the 
base manifold (also see section \ref{section4}) to construct a desired 
test function.  Since 
the normal bundle may not be trivial, we only consider the vertical 
functions that are radially symmetric on each fiber.

We state the  main results of this paper below.

\newtheorem{theorem}{Theorem}
\begin{theorem}\label{thm1}
Let $(F, f^*(ds_{E}^{2}))$ be an order-$k$ quantum tube with radius $r$,
with  the base
manifold $\varSigma$ satisfying condition A2).
Then the essential spectrum $\sigma_{ess}(\Delta)$ 
of the Dirichlet
Laplacian $\Delta$ can be estimated from below as 
\begin{equation} 
\inf\sigma_{ess}(\Delta) \geq \rho(k)^{2}/r^2,
\end{equation}
where $\rho(k)^2 > 0$ is the first eigenvalue of the Dirichlet 
Laplacian on the k-dimensional Euclidean ball of radius $r=1$.  
\end{theorem}    

\begin{theorem}\label{thm2}
Let  $(F, f^*(ds_{E}^{2}))$ be an order-$k$ 
quantum tube with radius r and base
manifold $\varSigma$ satisfying condition A1).  
We assume that $\varSigma$ is parabolic and not totally geodesic.  Suppose 
the function $\sum_{p=1}^{[n/2]} \mu_{2p}
K_{2p}$ is integrable and that  
\begin{equation}
\int_{\varSigma} \sum_{p=1}^{[n/2]} \mu_{2p}
K_{2p} \, \leq \, 0, 
\end{equation}  
where $\mu_{2p}$ are positive constant coefficients defined 
as $$\mu_{2p} =p(2p+k-2)
\int_0^1 t^{2p+k-3} 
\chi^2\, dt;$$
$\chi=\chi(r)$ defines the eigenfunction of the first eigenvalue of the Dirichlet Laplacian on the $k$-dimensional Euclidean ball of radius $1$:
\begin{equation}\label{psq}
\left\{
\begin{array}{l}
\chi''+\frac{n-1}{r}\chi=-\rho(k)^2\chi\\
\chi'(0)=0\\
\chi(1)=1
\end{array}
\right.;
\end{equation}
$K_{2p}$ are intrinsic curvature functions on $\varSigma$ defined
in Definition~\ref{def11}.   
Then the total spectrum 
$\sigma_0(\Delta)$ is strictly bounded above as
\begin{equation} 
\inf\sigma_0(\Delta) \, < \, \rho(k)^2/r^2. 
\end{equation}
\end{theorem}

\vskip 0.1cm

Combining the two theorems above we have the conclusive result:

\begin{theorem}\label{thm3}
Let $(F, f^*(ds_{E}^{2}))$ be an order-$k$ quantum tube with radius r 
and base
manifold $\varSigma$ satisfying conditions A1) and A2).  Moreover, 
we assume that $\varSigma$ is a parabolic manifold with 
$\sum_{p=1}^{[n/2]} \mu_{2p}K_{2p}$ integrable and 
\begin{equation*}
\int_{\varSigma} \sum_{p=1}^{[n/2]} \mu_{2p}
K_{2p} \, \leq \, 0, 
\end{equation*}
If $\varSigma$ is not totally geodesic, then the discrete spectrum of
 the Dirichlet Laplacian is non-empty.  
\end{theorem}

Of particular interest is the surface case. Corollary 1.1 of~\cite{ll-1} 
(or results in~\cite{cek-1})
can be generalized to higher codimension without any rephrasing.

\begin{corollary}\label{cor1}
Suppose that $\Sigma$ is a complete immersed surface of
$\R^n (n\geq 2)$ such that the second fundamental form $\vec A\rightarrow
0$. Suppose that the Gauss curvature is integrable and suppose
that 
\begin{equation}\label{1-4}
e(\Sigma)-\sum\lambda_i\leq 0,
\end{equation}
where $\lambda_i$ is the isoperimetric constant
at each end defined as follows. Let
$E_1,\cdots,E_s$ be the ends of the surface
$\Sigma$. For each  $E_i$, we define
\begin{equation}\label{iso}
\lambda_i=\underset{r\rightarrow\infty}{\lim}\,\frac
{A_i(r)}{\pi r^2},
\end{equation}
where $A_i(r)$ is the area of the ball $B(r)\cap
E_i$.
Let $a$ be a positive number such that $a||A||<C_0<1$.
If $\Sigma$ is not totally geodesic, then the ground
state of the quantum layer $\Omega$ exists. In particular, if
$e(\Sigma)\leq 0$, then the ground state exists.
\end{corollary}

Thus in the surface case the existence of discrete spectra can be controlled 
by the topology of the surface.  This fact and its possible generalizations to 
higher dimensions should be interesting to pursue further.

 {\bf Acknowledgement.} The authors would like to thank
 P. Exner,  
 D. Krej\v{c}i\v{r}\'{i}k, P. Li, and C-L. Terng for their interest and useful remarks.
 The first author
also likes to thank P.  Li for his financial support during the 
summer of 2004, and A. Klein for conveying much knowledge on functional analysis in relation to 
mathematical physics.  Much gratitude also goes to R. Mazzeo for many suggestions to  this paper, and for the 
information about the tube invariants of Weyl in Gray's book.

\section{Geometry of Quantum Tubes}\label{section2}

First we set up Fermi coordinate systems on quantum tubes to perform
local computations. The coordinates are defined in the book~\cite{gray}.

\begin{definition}
Consider the isometric immersion 
$\varSigma^{n} \hookrightarrow \mathbb{R}^{n+k}$ of the 
base manifold.
For a local coordinate chart $(U, \phi)$ on $\varSigma$, 
we trivialize 
the normal bundle $T^{\bot}\varSigma$ over $U$ with an orthonormal 
frame 
$\{\eta_1, \cdots, \eta_k\}$.  Then for each $x \in U$,  
$(x,\xi) \in T_{x}^{\bot}\varSigma$ we define local coordinates 
$(x_1, \cdots, x_n, u_1, \cdots, u_k)$, where
$\xi = u_{\alpha} \eta_{\alpha}(x)$.  We call such a coordinate system
a Femi coordinate system on $T^{\bot}\varSigma$.   \end{definition}

We now fix the convention of indexing the horizontal coordinates by 
Roman letters and fiber coordinates by Greek letters.  Then for
a Fermi coordinate system the entries in the metric tensor can be
expressed as
\begin{equation}
f^*(ds_{E}^{2}) = 
G_{ij}dx_{i}dx_{j} 
+ G_{i\alpha}dx_{i}du_{\alpha} 
+ G_{\alpha i}du_{\alpha}dx_{i} 
+ G_{\alpha\beta}du_{\alpha}du_{\beta}. 
\end{equation}
Then a straightforward but careful calculation gives 
\begin{eqnarray}\label{me}
\lefteqn{
\qquad  G_{ij} = 
g_{ij} - 2u_{\alpha}\langle S_{\eta_{\alpha}}(\partial_i), 
\partial_j\rangle} \nonumber\\ & & +  
u_{\alpha}u_{\beta}\langle 
S_{\eta_{\alpha}}(\partial_i), S_{\eta_{\beta}}(\partial_j)\rangle 
+  u_{\alpha}u_{\beta}\langle 
\nabla_{\partial i}^{\bot}\eta_{\alpha}, 
\nabla_{\partial j}^{\bot}\eta_{\beta}\rangle,
\end{eqnarray}
\begin{equation}\label{Gib}
G_{i\beta} = 
u_{\alpha}\langle \nabla_{\partial i}^{\bot}\eta_{\alpha}, 
\eta_{\beta}\rangle,
\end{equation}
\begin{equation}\label{I}
G_{\alpha\beta} = \delta_{\alpha\beta},
\end{equation}
where $\partial_i = \frac{\partial}{\partial x_i}$, and $g_{ij} dx_idx_j$ is the Riemann metric on
$\Sigma$. Note that by our notations
\begin{equation}\label{new-1}
u_{\alpha}u_{\beta}\langle 
\nabla_{\partial i}^{\bot}\eta_{\alpha}, 
\nabla_{\partial j}^{\bot}\eta_{\beta}\rangle
=G_{i\gamma}G_{j\gamma}.
\end{equation}

Compare to the metric tensor calculated in \cite{ll-1} 
for the quantum layer,
here we see the appearance of terms 
involving the normal connection.  Now, from the
Frobenius theorem (see Theorem 1.4.3 in \cite{pt}) we know 
that the flatness of the normal bundle (i.e.,$R^{\bot} 
\equiv 0$ on $\varSigma$) guarantees that we can choose normal frames on 
$\varSigma$ to annihilate these terms.  The existence of these parallel 
normal frames (although not necessarily orthonormal) would 
simplify our analysis as it would reduce the metric 
tensor to block form similar to that of in \cite{ll-1}.  However, 
as remarked in the introduction, this is too much to ask of immersions.    

\begin{remark}
Formulas~\eqref{me} through~\eqref{new-1} are also valid for the case 
where $\varSigma$ is a curve ($n=1$).  It is interesting to note that 
in the case of a unit speed curve in $\mathbb{R}^3$, where we employ the 
usual Frenet frames, the terms $\nabla_{\partial i}^{\bot}\eta_{\alpha}$  
contain the torsion of the curve.  One may 
wish to carry this intuition (and intuition only) further 
in the cases of our concern, namely 
$n\geq 2$.  Then we will eventually see that the presence of torsion in the 
metric is in fact irrelevant for the existence of discrete spectra.  Now, 
there is a series of pre-existing work on the spectral analysis of tubes 
built over curves (cf. ~\cite{es-1,gj,de-2,cdfk,kr-1}).  A notable paper is 
~\cite{cdfk}, where the authors essentially 
demonstrated that the only criterion for the existence of discrete spectra 
is for the tube to be not straight but asymptotically straight.  There they 
also proved that the essential spectrum is identical to that of the straight 
tube.  It seems that the more intrinsic methods in this paper 
can also be used to analyze the curve case, but perhaps it is better to 
discuss this in a future paper.        
\end{remark}

From~\eqref{me} we see that the $n\times n$ matrix $(G_{ij})$ can be 
written as
\begin{equation}\label{Gij}
(G_{ij}) = (\tilde{G}_{ij}) + (\mathcal{N}_{ij}),
\end{equation}
where 
\begin{equation}\label{tildeG}
(\tilde{G}_{ij}) = (I - u_{\alpha}H_{\alpha})^{2}(g_{kj}),
\end{equation}
and
\begin{equation}\label{Nij}
(\mathcal{N}_{ij}) = (u_{\alpha}u_{\beta}\langle 
\nabla_{\partial_i}^{\bot}\eta_{\alpha}, 
\nabla_{\partial_j}^{\bot}\eta_{\beta}\rangle)
\end{equation}
for the induced metric tensor $(g_{ij})$ on $\varSigma$ and the 
local matrix $H_{\alpha}$ representing the shape operator 
$S_{\eta_{\alpha}}$\footnote{That is, $(H_\alpha)_{ij}=\langle S_{\eta_\alpha}(\partial_i),\partial_k\rangle g^{kj}$.}.

Next, we let $(\tilde{G}_{ij})_{n\times n} = \tilde{G}$ and 
$(G_{i\beta})_{n\times k} = C$.  Then in view of (\ref{I}), ~\eqref{new-1},
~(\ref{Gij}), and (\ref{tildeG}) we see that the metric tensor $G$
on $T^{\bot}\varSigma$ with respect to the Fermi coordinate
system has the block form
\begin{equation}\label{meT}
\left(
\begin{matrix}
\tilde{G} + CC^{T} & C\\
C^{T} & I
\end{matrix}
\right).
\end{equation}
The matrix $C$ contains the normal connection terms, but we
will see that the resulting volume element of the quantum tube
will not depend on the normal connection explicitly.  

\vskip 0.1cm

\newtheorem{lemma}{Lemma}
\begin{lemma}\label{lem1} Using the notations as above, we have
$\det G = \det\tilde{G}$
\end{lemma}

\text{\bf{Proof.}} 
The trick is the following multiplication of matrices
\begin{equation}\label{inte}
\left(
\begin{matrix}
\tilde{G} +CC^T & C \\
C^T & I
\end{matrix}
\right)
\left(
\begin{matrix}
  I & 0 \\
  -C^T & I
\end{matrix}
\right) = 
\left(
\begin{matrix}
  \tilde{G} & C \\
  0 & I 
\end{matrix}
\right).
\end{equation}

\qed

\begin{corollary}
The condition $r \leq (\sqrt{k} \hskip 0.05cm \epsilon_{o})^{-1}$ 
guarantees that the map $f$ given by (1) is an immersion when 
restricted
to the sub-bundle $F$, where $(F,f^*(ds_{E}^{2}))$ is an order-k 
quantum
tube with radius r.  
\end{corollary}

\text{\bf{Proof.}} 
We know that the requirement that $df$ be nonsingular at a point in 
$T^{\bot}\varSigma$ is equivalent to the condition that the pullback
metric tensor $G$ be a nonsingular matrix at the same point. By Lemma
\ref{lem1}, this is reduced to knowing the invertibility of the matrix
$\tilde{G}$.  Suppose $\tilde{G}$ is invertible at a point, then

\begin{equation}
\tilde{G}^{-1} = (I - u_{\alpha}H_{\alpha})^{-2} g^{-1}
\end{equation} 
               
from (\ref{tildeG}), where $g$ denotes $(g_{ij})$.  Note that
$I - u_{\alpha}H_{\alpha}$ is invertible because $g$ is.  
Now, we can write the Taylor expansion
\begin{equation}
(I - u_{\alpha}H_{\alpha})^{-1} = 
\sum_{n=0}^{\infty}u_{\alpha}^{n}(H_{\alpha})^{n},
\end{equation}
which converges if and only if 
\begin{equation}\label{suff}
\|\sum^k_{\alpha = 1} u_{\alpha}H_{\alpha}\| < 1.
\end{equation}
On the other hand, we observe that 
\begin{equation*}
\|\sum^k_{\alpha = 1} u_{\alpha}H_{\alpha}\| \hskip 0.3cm \leq 
\hskip 0.3cm  
\sum^k_{\alpha = 1} |u_{\alpha}|\|H_{\alpha}\| \hskip 0.3cm       
\leq \hskip 0.3cm 
\|\vec{A}\|\sum^k_{\alpha = 1}|u_{\alpha}| \hskip 0.3cm 
\leq \hskip 0.3cm 
\epsilon_{o}\sum^k_{\alpha = 1}|u_{\alpha}|
\end{equation*}
by condition A1) and the relation between the shape
operator and the second fundamental form.

Therefore a sufficient condition for guaranteeing (\ref{suff}) is 
$\sum^k_{\alpha = 1}|u_{\alpha}| < (\epsilon_{o})^{-1}$.  
However, to 
yield a geometrically invariant condition we further 
observe by the Cauchy inequality that
\begin{equation*}
\sum_{\alpha = 1}^k |u_{\alpha}| \hskip 0.3 cm
 \leq \hskip 0.3 cm 
\sqrt{k} \hskip 0.05cm 
\left(
\sum_{\alpha = 1}^k 
|u_{\alpha}|^2
\right)^{\frac{1}{2}}.
\end{equation*}
Thus we also obtain the sufficient condition $\left(
\sum_{\alpha = 1}^k 
|u_{\alpha}|^2
\right)^{\frac{1}{2}} < 
(\sqrt{k}$ \hskip 0.05cm $\epsilon_{o})^{-1}$ that guarantees 
(\ref{suff}), 
and this is a statement about the distance to $\varSigma$, which
is invariant under the choice of orthonormal frame in the normal
bundle.  Then by the definition of $r$, the condition
$r \leq (\sqrt{k} \hskip 0.05cm \epsilon_{o})^{-1}$ ensures that 
the map $f$ given by (\ref{E:imm}) is an immersion when restricted to 
the 
quantum tube. 

\qed

\begin{corollary}
Assume that $r<\frac{1}{\sqrt k}$. 
Then the  volume element of an order-$k$ quantum tube satisfies
\begin{equation} \label{E:cor2}
 (1 - \epsilon_{o})^k\, du_{A} \,d\varSigma 
 \hskip 0.2cm 
 \leq \hskip 0.2cm 
 \det (I-u_{\alpha}H_{\alpha})\, du_{A} \,d\varSigma
 \hskip 0.2cm  \leq 
 \hskip 0.2cm
(1 + \epsilon_{o})^k\, du_{A} \,d\varSigma,
\end{equation}
where $\epsilon_{o}$ is the constant that appears in condition A1),
and we assume $\epsilon_0<1$.
\end{corollary}

\text{\bf{Proof.}}
By Lemma \ref{lem1} and equation (\ref{tildeG}) we see that 
\begin{eqnarray*}
  \lefteqn{
  dF = \det(I - u_{\alpha}H_{\alpha}) \sqrt{g}du_A\,dx_{J}} \\
  & & = \det(I - u_{\alpha}H_{\alpha}) du_A \, d\varSigma.
\end{eqnarray*}
The corollary follows from the fact that $\|u_\alpha H_\alpha\|<\epsilon_0$.
\qed

We remark here that if the norm of the 
second fundamental form is bounded 
by any constant $0 < \epsilon /\sqrt k<\epsilon_0< 1$, (\ref{E:cor2}) still holds and we 
have 
\begin{equation}\label{E:ncor2}
(1 - \epsilon)^k\, du_{A}\, d\varSigma 
\hskip 0.2cm 
\leq \hskip 0.2cm 
\det (I-u_{\alpha}H_{\alpha})\, du_{A}\, d\varSigma \hskip 0.2cm  \leq 
\hskip 0.2cm
(1 + \epsilon)^k\, du_{A}\, d\varSigma.
\end{equation}
The inequalities above will be key in the estimate of the lower 
bound of the essential spectrum.  

We also need  to calculate the inverse matrix $G^{-1}$.  A straightforward computation using~\eqref{inte} gives
\begin{align}\label{G-1}
\begin{split}
G^{-1} &= 
\left(
\begin{matrix}
  I & 0 \\
  -C^T & I
\end{matrix}
\right)
\left(
\begin{matrix}
\tilde{G}^{-1} & -\tilde{G}^{-1}C \\
0 & I 
\end{matrix}
\right)\\&= 
\left(
\begin{matrix}
  \tilde{G}^{-1} & -\tilde{G}^{-1}C \\
  -C^{T}\tilde{G}^{-1} & C^{T}\tilde{G}^{-1}C + I
\end{matrix}
\right)\\&= 
\left(
\begin{matrix}
\tilde{G}^{-1} & -\tilde{G}^{-1}C \\
  -C^{T}\tilde{G}^{-1} & C^{T}\tilde{G}^{-1}C 
\end{matrix}
\right) + 
\left(
\begin{matrix}
0 & 0 \\
0 & I
\end{matrix}
\right).
\end{split}
\end{align}

\qed
 
We will end this section with one property of the above equation:

\vskip 1cm 

\begin{lemma}\label{lem2}
$\left(
\begin{matrix}
\tilde{G}^{-1} & -\tilde{G}^{-1}C \\
  -C^{T}\tilde{G}^{-1} & C^{T}\tilde{G}^{-1}C 
\end{matrix}
\right)$ is positive semidefinite.
\end{lemma}

\text{\bf{Proof.}}\hskip 0.3cm
Let $\vec{v} = (v_1, v_2,\cdots, v_n, v_{n+1},\cdots,v_{n+k}) = 
(\vec{v}_i \hskip 0.3cm \vec{v}_{\alpha}) 
\in \mathbb{R}^{n+k}$, where $\vec{v}_i$ 
denotes the first $n$ components
and $\vec{v}_{\alpha}$ denotes the last $k$ components.  
Then we see that
\begin{eqnarray}\label{psem}
\lefteqn{
\left(
\begin{matrix}
\vec{v}_i & \vec{v}_{\alpha}
\end{matrix}
\right)
 \left(
 \begin{matrix}
 \tilde{G}^{-1} & -\tilde{G}^{-1}C \\
  -C^{T}\tilde{G}^{-1} & C^{T}\tilde{G}^{-1}C 
 \end{matrix}
 \right)
\left(
\begin{matrix}
\vec{v}_i \\
\vec{v}_{\alpha}
\end{matrix}
\right)} \nonumber \\
 & & = 
\tilde{G}^{ij}v_{i}v_{j} - 
2\tilde{G}^{ij}G_{j\alpha}v_{i}v_{\alpha}
 + G_{\alpha i}\tilde{G}^{ij}G_{j\beta}v_{\alpha}v_{\beta} 
\nonumber \\
 & & = 
(\vec{v}_{i} - C\vec{v}_{\alpha})^{T}\tilde{G}^{-1}
(\vec{v}_{i} - C\vec{v}_{\alpha}).
\end{eqnarray}
Since the matrix  $\tilde G$ is clearly positive definite, the lemma is proved.

\qed

\section{Lower Bound Estimate of the Essential Spectrum}\label{section3}
Let $f\in C_0^\infty(F)$. In the Fermi coordinate system, we have
\begin{equation}\label{fib}
|\nabla f|^{2} \geq |\nabla^{\bot}f|^2 
\end{equation}
by Lemma \ref{lem2}.  The fiber-wise gradient $\nabla^{\bot}$ is 
essentially the gradient on the $k$-dimensional Euclidean ball. 

It is clear that the boundary of a quantum tube $F$
with radius 
$r$ and base manifold $\varSigma$ has as boundary the smooth manifold
$\varSigma \times \partial B^{k}(0,r)$, where we identify 
the $k$-dimensional Euclidean open ball 
$B^k(0,r)$ with the fiber 
$F(x)$ at each $x \in \varSigma$.  Then since $f = 0$ on 
$\partial F$, we have the Poincar\'e inequality
\begin{equation}\label{Poinc}
\int_{B^k(0,r)} |\nabla^{\bot}f|^2 \hskip 0.1cm du_{A} 
\hskip 0.3cm \geq \hskip 0.2cm  
\frac{\rho(k)^2}{r^2}\int_{B^k(0,r)}f^2 \hskip 0.1cm du_{A},
\end{equation} 
where $\rho(k)^2$ is defined in Theorem~\ref{thm1}.

Now we estimate the essential spectrum form below and prove
Theorem \ref{thm1}.

\text{\bf{Proof of Theorem~\ref{thm1}.}} \hskip 0.3cm 
Fix a reference point $x_{o}$ on $\varSigma$, condition A2)
implies that for any 
$0 < \epsilon < 1$ there exists an open ball 
$B(x_{o}, R) \subset \varSigma$ 
such that (\ref{E:ncor2}) holds outside its closure.  
We may then consider the compact set
\begin{equation*} 
K = \overline{B(x_{o}, R)} \times \overline{B^{k}(0, t)}
\end{equation*}
in $F$, where  
$t < r$.  Using (\ref{E:ncor2}) and (\ref{fib}) we obtain

\begin{equation}\label{Dirineq}
\frac{\int_F |\nabla f|^2}{\int_F f^2} \geq
\left(\frac{1 - \epsilon}{1 + \epsilon}\right)^k
\left(
\frac{\int_F |\nabla^{\bot}f|^2 \hskip 0.2cm du_{A}\, d\varSigma}
{\int_F f^2 \hskip 0.2cm du_{A}\, d\varSigma}\right)
\end{equation}

for every $f \in C_{0}^{\infty}(F\setminus K)$. 
Then (\ref{Poinc}) implies that 
\begin{eqnarray}
\lefteqn{
\int_F |\nabla^{\bot} f|^2 du_{A}\, d\varSigma = 
\int_{\varSigma}\int_{B^k(0,r)} 
|\nabla^{\bot} f|^2 du_{A}\, d\varSigma} \nonumber \\  
 && \hskip 2.1cm \geq \hskip 0.2cm  
\frac{\rho(k)^2}{r^2} \int_{\varSigma}\int_{B^k(0,r)} f^2 \hskip 0.1cm 
du_{A}\, d\varSigma,   
\end{eqnarray}         
from which (\ref{Dirineq}) becomes
\begin{equation}
\frac{\int_F |\nabla f|^2}{\int_F f^2} > 
\left(\frac{1 - \epsilon}{1 + \epsilon}\right)^k\frac{\rho(k)^2}{r^2}.
\end{equation}
for every $f \in C_{0}^{\infty}(F\setminus K)$.  Then in view of  
(\ref{infess}) we obtain
\begin{equation}
\inf\sigma_{ess}(\Delta) \geq 
\left(\frac{1 - \epsilon}{1 + \epsilon}\right)^k\frac{\rho(k)^2}{r^2}.
\end{equation}

Since $\epsilon$ is arbitrarily small, we obtain the desired estimate.  

\qed

This result generalizes the lower bound estimate of the essential 
spectrum obtained in the codimension-$1$ case in \cite{ll-1}, where 
$\rho(1)^2 /r^2= \kappa_{1}^{2} = \frac{\pi^2}{4r^2}$.

\vskip 1cm

\section{Parabolicity of Complete Manifolds}\label{section4}
In this section $M$ will denote a complete manifold without boundary.  
The {\it parabolicity} on $M$ 
was first defined as an analytic notion related to PDEs.  Recall that a 
Green's function is the fundamental solution to Poisson's equation.   

\begin{definition}
M is parabolic if it does not admit any positive Green's function, otherwise 
it is said to be nonparabolic.
\end{definition}

The definition above and the subsequent material in this section is based 
on the  survey paper by P.~Li (\cite{li}).  

The definition given above for parabolicity may seem quite intangible.  
However, P.~Li and L.~F.~Tam in \cite{li-tam} gave a procedure for 
constructing a Green's function on $M$.  From the procedure, 
which involves applying the maximum principle to the 
positive Green's function  
on each $\Omega_i$ of a compact exhaustion $\{\Omega_i\}$ of $M$ and 
analyzing the limit of such sequence of  Green's functions, they could 
extract the following equivalent condition for parabolicity 
(see \cite{li}):

\newtheorem{proposition}{Proposition}
\begin{proposition}\label{prop1}
Let $B(s)$ be a geodesic ball of radius $s$ in $M$ centered at any fixed 
point $x_0$.  Let $R>s>1$.  Then let $\psi_R$ be the solution to 
the following problem 
\begin{equation*}
\begin{cases}
 \Delta \psi = 0 \quad & \text{on $B(R)\setminus B(s)$};\\
 \psi|_{B(s)} \equiv 1; \\
  \psi|_{\varSigma \setminus B(R)} \equiv 0.
\end{cases}
\end{equation*}
\end{proposition}  
Then $M$ is parabolic if and only if 
\begin{equation}\label{gotozero}
\int_{M}|\nabla\psi_R|^2 \rightarrow 0 \quad 
\text{as $R \rightarrow \infty$}.
\end{equation}

The harmonic functions $\{\psi_R\}$ is the key use of parabolicity 
in this paper, 
and will be used directly in the next section.  The characterization of 
parabolicity is still anaytical in the proposition above.  However, very 
often $M$ being parabolic is also a geometric condition.  First of all, there 
is the following result proven independently by Grigor'yan 
\cite{gro1,gro2} and Varopoulos \cite{v2}:

\begin{theorem}\label{var}   
Let $V(t)$ be the volume of the geodesic ball $B(t)$ centered 
at any $p\in M$.  If $M$ is nonparabolic, then  
\begin{equation}\label{volgrow}
\int_1^{\infty}\frac{tdt}{V(t)} < \infty.
\end{equation}
Therefore, if a geodesic ball is at most of quadratic growth then 
$M$ must be parabolic.  
\end{theorem}

In particular, the following corollary holds:

\begin{corollary}
Any smooth surface with integrable Gaussian curvature must be parabolic.
\end{corollary}

The corollary above is the reason why we did not have to 
assume the base manifold is parabolic in Corollary \ref{cor1}.  

While the converse to Theorem~\ref{var} is not true by a counter-example 
of Greene (see \cite{v2}, or more directly \cite{li}), there is for example 
the following result.

\begin{theorem}[Varopoulos \cite{v1}]
If $M$ has non-negative Ricci 
curvature, then M is non-parabolic if and 
only if the volume growth condition~(\ref{volgrow}) holds.
\end{theorem}    

Although its use is unnecessary, the theorem above implies that 
with the standard Euclidean metric, $\mathbb{R}^2$ is parabolic 
while $\mathbb{R}^n$ is nonparabolic  for all $n \geq 3$.  The reader 
can also view this distinction between $\mathbb{R}^2$ and 
$\mathbb{R}^n$($n \geq 3$) as a motivation for parabolicity (see 
our previous paper \cite{ll-1}).      In gereral, if the volume growth of a complete manifold is quadratic, then it is parabolic.

\section{Upper Bound Estimate of the Lower Bound of the Total Spectrum}\label{section5}
We already saw in the previous two sections that the appearance of 
the normal connection terms $G_{i\alpha}$ did not stop us from 
generalizing the corresponding estimate for the codimension-$1$ case.  
The same will be true for the total spectrum.  

First we recall from the introduction  

\begin{equation}\label{inftot}
\inf\sigma(\Delta) = 
\inf_{f \in C_{0}^{\infty}(F)}
\frac{\int_F |\nabla f|^2}{\int_F f^2}.
\end{equation}  

By rescaling, we can assume the radius of $F$ is $1$.
To prove Theorem \ref{thm2}, it then suffices to find a test function
$\phi$ smooth almost everywhere on $F$ such that 

\begin{equation}\label{opq}
\int_F |\nabla\phi|^2 - \rho(k)^2\int_F \phi^2 
\hskip 0.2cm < \hskip 0.2cm 0.
\end{equation}
As in \cite{ll-1}, the prototype test function is of  
the form $\phi = \chi\psi$, where $\chi$ depends 
only on the fiber and $\psi$ depends only on the base manifold.  
Then since 
$\nabla\chi\psi = \chi\nabla\psi + \psi\nabla\chi$, we see that 

\begin{equation}
|\nabla\phi|^2 = |\nabla\chi\psi|^2 = 
\chi^2|\nabla\psi|^2 + \psi^2|\nabla\chi|^2+2\psi\chi
\langle\nabla\psi,\nabla\chi\rangle.
\end{equation}

Then (\ref{opq}) becomes 
\begin{equation}\label{totineq}
\int_F \psi^2(|\nabla\chi|^2 - \rho(k)^2\chi^2) + 
2\int_F\psi\chi\langle\nabla\psi,\nabla\chi\rangle
+
\int_F \chi^2|\nabla\psi|^2 \hskip 0.2cm < 
\hskip 0.2cm 0.
\end{equation}
           
We will assume that the base manifold is parabolic, so for a fixed point 
$x_{o} \in \varSigma$ and any 
$R > s >1$, we can let the horizontal function be $\psi = \psi_R$ in 
Proposition~\ref{prop1}, satisfying~(\ref{gotozero}).

The choice of the fiber-wise function $\chi$ will also be similar to 
those in \cite{ll-1}.  As we said in the introduction, 
$\chi$ should be radially symmetric in the fibers.  
Let $t$ be the parameter for the length of each 
vector $\eta \in T^{\bot}_{x}\varSigma$.
Since locally each fiber $F(x)$ is identified with $B^k(0,1)$, 
We will assume $\chi$ is of the form $\chi (t)$ on $B^k(0,1)$. 
  
Then a straightforward computation using ~\eqref{G-1} shows that
\begin{equation}\label{delchi}
|\nabla\chi|^2 = 
\left(
\frac{u_{\alpha}u_{\beta}}{t^2}G_{\alpha i}\tilde{G}^{ij}G_{j\beta} 
 \ + \  1
\right) 
\left| \frac{d\chi}{\partial t}\right| ^2.
\end{equation}
However, paying particular attention to the repeated indices 
for summation and using (\ref{Gib}), we see that for each 
$j = 1,2,\cdots, n$,
\begin{eqnarray}\label{wonderzero}
\lefteqn{
u_{\beta} G_{j\beta} = u_{\beta}u_{\gamma}
\langle \nabla_{\partial_j}^{\bot}\eta_{\gamma}, 
\eta_{\beta}\rangle}\nonumber \\ 
& & \hskip 0.4cm =
\sum_{\beta\ne\gamma}u_{\beta}u_{\gamma}
\langle \nabla_{\partial_j}^{\bot}\eta_{\gamma}, 
\eta_{\beta}\rangle\nonumber \\ 
& & \hskip 0.4cm = 0 ,
\end{eqnarray}
since $\{\eta_1, \cdots, \eta_k\}$ is an orthonormal frame and 
symmetry implies $\langle\nabla_{\partial_j}^{\bot}\eta_{\gamma}, 
\eta_{\beta}\rangle = -\langle\nabla_{\partial_j}^{\bot}\eta_{\beta}, 
\eta_{\gamma}\rangle$. 
Therefore (\ref{delchi}) becomes
\begin{equation}\label{ndelchi}
|\nabla\chi|^2 = \left| \frac{d\chi}{dt}\right|^2.
\end{equation} 
Using the same reasoning, we have
\begin{equation}\label{ndelchi-1}
\langle\nabla\psi, \nabla\chi\rangle = 0
\end{equation}
if $\chi$ is rotationally symmetric.

Therefore, equation~(\ref{wonderzero}) again shows that the 
normal connection plays no explicit role in the analysis.  
Equation~(\ref{totineq}) now takes 
the familiar 
form 
as in \cite{ll-1}: 
\begin{equation}\label{ntotineq}
\int_F \psi^2
\left(
\left|
\frac{d\chi}{dt}
\right|^2 
- \rho(k)^2 \chi^2
\right) 
+ \int_F \chi^2 |\nabla\psi|^2 
\hskip 0.1cm < \hskip 0.1cm 0.
\end{equation} 

For the purpose of this section, we will be integrating in 
polar coordinates on the fibers.  The use of polar coordinates 
is of course natural since our fibers are essentially Euclidean 
balls, although such a use 
was not necessary in proving Theorem 1.  Now recall that for a point 
$(x,\xi)\in F$, we have
\begin{equation}\label{detIH}
\det (I - u_{\alpha}H_{\alpha})  = \,   
\sum_{j=0}^{n}(-1)^j  C_j(S_{\xi}).
\end{equation}
In view of this we can write $\xi = t\eta$, where
$\eta = \xi/\|\xi\|$.  Then $S_{\xi} = tS_{\eta}$ and 
(\ref{detIH}) becomes
\begin{equation}\label{ndetIH}
\det (I - u_{\alpha}H_{\alpha}) \, =
\, \det (I - tH_{\eta}) \, = \,   
\sum_{j=0}^{n}(-1)^j t^j C_j(S_{\eta}),
\end{equation}
where $H_{\eta}$ is the matrix of $S_{\eta}$ with respect 
to the Fermi coordinate system and 
we define $C_0(S_\eta) = 1$.

The quantity $C_{j}(S_{\eta})$ is a function on the unit sphere
bundle
\begin{equation*} 
F_1 = \{(x,\eta) \in T^{\bot}\varSigma \mid \|\eta\| = 1\} 
\subset F.
\end{equation*}
\newtheorem{Definition}{definition}
\begin{definition}\label{def11} \hskip 0.2cm
We define the following function on $\varSigma$ 
(over all $\eta \in F_1(x)$):
\begin{equation*}
K_j = \int_{S^{k-1}}C_j(S_\eta)d\sigma,
\end{equation*}
and call it the jth-curvature of the quantum tube $F$,  
where $S^{k-1}$ is identified with the fiber $F_1(x)$ and 
the integration is with respect to the boundary measure 
$d\sigma$ induced by the orientation of $F$.
\end{definition}
An immediate 
observation from the definition is that all the 
odd-curvatures of $F$ are 0.  This is because we have 
$C_j(S_{-\eta}) = -C_j(S_\eta)$ 
and $S^{k-1}$ is radially 
symmetric.  

Let $\cal R$ be the curvature operator of $\Sigma$. 
We define
the $p$-th trace  of  $\cal R$ as
\begin{align}\label{trRp}
\begin{split}
&\qquad tr(\mathcal{R}^{p}) = 
\frac{1}{2^p ((2p)!)^2}\\
&\times
\sum_{I} 
\sum_{\sigma, \tau \in S_{2p}} sgn(\sigma\tau) 
R_{i_{\sigma(1)}i_{\sigma(2)}i_{\tau(1)} i_{\tau(2)}} \cdot\cdot\cdot 
R_{i_{\sigma(2p-1)}i_{\sigma(2p)}i_{\tau(2p-1)}i_{\tau(2p)}}.
\end{split}
\end{align} 

Then we have the following result from~\cite[equation (4.15)]{gray}:

\begin{proposition}  
For an isometric immersion of a manifold 
$\varSigma$ with any codimension $k$ 
into Euclidean space, we have at each point on $\varSigma$  
\begin{equation}\label{C2case}
K_{2p}=\int_{S^{k-1}}C_{2p}(S_{\eta}) = 
\frac{(2p)!\pi^{k/2}}{2^{2p-1}p!\Gamma(p+k/2)}\, tr(\mathcal{R}^p).
\end{equation}
\end{proposition}

\qed

The formula above is nice in that up to a multiple of a constant 
depending only on the dimesion and codimension, $K_{2p}$ is determined 
entirely by the intrinsic Riemmanian structure of $\varSigma$.

\vskip 0.5cm

{\bf Proof of Theorem~\ref{thm2}.}
Using (\ref{me}), condition A1), and (\ref{G-1}), 
we can show that  
\begin{equation}\label{hestm}
\int_F \chi^2|\nabla\psi|^2 \, \leq \, 
C_{1}\int_{\varSigma} |\nabla_{\varSigma}\psi|^2
\end{equation}
for some constant $C_1 > 0$.

We consider the case 
\begin{equation}\label{<0}
\int_{\varSigma}\sum_{p=1}^{[n/2]}\mu_{2p}K_{2p}\, d\varSigma
\, < \, 0.
\end{equation}
first.  Here we pick the prototype test function $\chi\psi$ 
described earlier. 

Let $\rho$ denote $\rho(k)$ and 
$B(0,1)$ denote the unit open ball in $\mathbb{R}^k$.   
Then  using (\ref{ndetIH}) we have
\begin{eqnarray}\label{major}
\lefteqn{
\int_F \psi^2 
\left(
\left|
\frac{d\chi}{dt} 
\right|^2 - \rho^2\chi^2
\right)
det(I - u_{\alpha}H_{\alpha}) \, du_{A}\, d\varSigma} \nonumber \\
& & = \int_\varSigma \psi^2 
      \int_{B(0,1)}  
\left(
\left|
\frac{d\chi}{dt} 
\right|^2 - \rho^2\chi^2
\right)
\sum_{j=0}^{n}(-1)^{j}t^{j}C_j(S_\eta)\, du_{A}\, d\varSigma
\nonumber \\ 
& & = 
\int_{\varSigma}\psi^2
\sum_{j=1}^{n} (-1)^j 
\int_{0}^{1}\int_{\partial B(0,t)}
t^j
\left(
\left|
\frac{d\chi}{dt} 
\right|^2 - \rho^2\chi^2
\right)
C_j(S_\eta)\, d\sigma_t \, dt \, d\varSigma \nonumber \\
& & = 
\int_{\varSigma}\psi^2
\sum_{j=1}^{n}(-1)^j
\int_{0}^{1}
t^j
\left(
\left|
\frac{d\chi}{dt} 
\right|^2 - \rho^2\chi^2
\right) dt
\int_{\partial B(0,t)} C_j(S_{\eta})\, 
d\sigma_{t}  \, \, d\varSigma \nonumber \\ 
& & = 
\int_{\varSigma}\psi^2
\sum_{j=1}^n (-1)^j \int_0^1
t^{j+k-1}
\left(
\left|
\frac{d\chi}{dt} 
\right|^2 - \rho^2\chi^2
\right) \, dt 
\int_{S^{k-1}}C_j(S_\eta)\, d\sigma \, \, d\varSigma
\nonumber \\
& & = 
\int_{\varSigma}\psi^2\sum_{p=1}^{[n/2]}\mu_{2p}
K_{2p}\, d\varSigma, 
\end{eqnarray} 
where for the fourth equality we used the fact that
\begin{equation}
\int_{\partial B(0,t)} f \,  d\sigma_t \, = \, 
t^{k-1}\int_{S^{k-1}} f \, d\sigma,
\end{equation}
if $f$ is independent of the radius parameter; and for the fifth equation
we used ~\eqref{psq}.

By our choice of $\psi$, (\ref{<0}), and in view of (\ref{major}), 
we can choose $s$ large enough so that
\begin{equation}\label{estm1}
\int_F \psi^2 \left(
\left|
\frac{d\chi}{dt}
\right|^2 
- \rho^2\chi^2
\right) \, \, 
< \, \, -\delta.
\end{equation} 
Moreover, by (\ref{hestm}) and (\ref{gotozero}), we can then
choose $R$ large enough so that 
\begin{equation}\label{estm2}
\int_F \chi^2|\nabla\psi|^2 \, \, < \, \, 
\delta.
\end{equation}
Combining (\ref{estm1}) and (\ref{estm2}) above yields 
(\ref{ntotineq}), which
proves the first part of the theorem.

Next we consider the case 
\begin{equation}\label{=0}
 \int_{\varSigma}\sum_{p=1}^{[n/2]}\mu_{2p}K_{2p}\, d\varSigma
 \, = \, 0.
\end{equation}
For brevity we will define the quadratic form $Q$ by
\begin{equation}\label{quadform}
 Q(f,g) = \int_F \langle \nabla f, \nabla g \rangle 
          \, - \, \rho^2 \int_F fg.
\end{equation}
Then we consider as in \cite{ll-1} (and \cite{dek}) the same test 
function
\begin{equation}
 \phi_{\epsilon} = \phi + \epsilon j\chi_1,
\end{equation}
where $\epsilon$ is a small number,  
$\phi = \chi\psi$ is the prototype test function, 
$j$ a smooth function supported in $B(s) \subset \varSigma$, and 
$\chi_1$ is some smooth function on $B(0,1)$ that vanishes
on the boundary.  
In view of (\ref{inftot}) and (\ref{quadform}), it suffices to show 
that
$Q(\phi_{\epsilon}, \phi_{\epsilon}) \, < \, 0$.

A quick computation shows that
\begin{equation*}
 Q(\phi_{\epsilon}, \phi_{\epsilon}) \, = \, 
 Q(\phi,\phi) + 2\epsilon Q(\phi,j\chi_1) + 
 \epsilon^2Q(j\chi_1,j\chi_1).
\end{equation*}

As in the previous case, we still have   
\begin{equation}\label{estmQphi}
 Q(\phi,\phi) \, \leq \, 
 \mathit{C}_1\int_{\varSigma}|\nabla\psi|^2\ + 
 \mathit{C}_2\int_{\varSigma}\psi^2 \sum_{p=1}^{[n/2]}\mu_{2p}K_{2p}.
\end{equation}
Now, by the fact that supp($j$) $\subset$ $B(s)$, and a 
careful calculation using the metric of $F$ given by (\ref{meT}), we
see that
\begin{equation}\label{23ndQ}
 Q(\phi, j\chi_1) = \int_\varSigma j\int_{B(0,1)} 
 \left(
 \frac{\partial\chi}{\partial u_{\alpha}}
 \frac{\partial\chi_1}{\partial u_{\alpha}} - 
 \rho^2\chi\chi_1
 \right) \det(I - u_{\alpha}H_{\alpha})\, 
 du_{\alpha}\, d\varSigma.
\end{equation}
Substituting the identity 
 $\sum_{\alpha = 1}^{k}
 \frac{\partial^2\chi}{\partial u_{\alpha}^2} = -\rho^2\chi$ 
on $B(0,1)$ into (\ref{23ndQ}), and integrating by parts we see that 
\begin{eqnarray}\label{2ndQ}
\lefteqn{
 Q(\phi, j\chi_1) = 
 -\int_{\varSigma}j\int_{B(0,1)} \chi_1
 \frac{\partial\chi}{\partial u_{\alpha}}
 \frac{\partial}{\partial u_{\alpha}}
 \det(I - u_{\alpha}H_{\alpha})\, du_{\alpha}\, d\varSigma}
 \nonumber \\
 & & \hskip 0.8cm = 
 -\int_{\varSigma}j\int_{B(0,1)}\chi_1
 \left< \nabla^{\bot}\chi , 
 \nabla^{\bot}\det (I - u_{\alpha}H_{\alpha}) \right> 
 \, du_{\alpha}\, d\varSigma. 
\end{eqnarray}
Note that we used the fact that $\chi$ being radially symmetric 
implies it must have derivatives $0$ at the origin.  

From (\ref{ndetIH}) we see that 
$\det (I - u_{\alpha}H_{\alpha}) = \det (I - tH_{\eta})$, 
which is a polynomial in $t$ at each $x \in \varSigma$.  
The assumption that $\varSigma$ is not totally geodesic implies 
that there exist a point $x \in \varSigma$ on which there is an
$\eta$ such that $H_{\eta} \ne 0$.  Then we can deduce that there
must be a $t_0 \in (0,r)$ such that 
$\frac{\partial}{\partial t}\det (I - tH_{\eta})(t_0) \ne 0$ 
at $x$.  
Moreover, since the zero set of $\chi$ is discrete, we can choose 
$t_0$ such that $\frac{d\chi}{dt}(t_0) \ne 0$ as well.

Now, using polar coordinates on $B(0,r)$ we have
\begin{equation*}
 \left< \nabla^{\bot}\chi , 
 \nabla^{\bot}\det (I - u_{\alpha}H_{\alpha}) \right> = 
 \frac{d\chi}{dt}\frac{\partial}{\partial t}\det (I - tH_{\eta}), 
\end{equation*}
which at $(x, t_0\eta)$ is not zero.  Then we simply choose 
$\chi_1$ to be a bump function about $t_0\eta$ 
and $j$ a bump function about $x$, 
so that $Q(\phi, j\chi_1) \ne 0$.  
Then we may choose a negative or positive $\epsilon$ small enough 
so that 
\begin{equation}\label{1stestm}
2\epsilon Q(\phi, j\chi_1) + \epsilon^2 Q(j\chi_1, j\chi_1) 
\, < \, -\delta.           
\end{equation} 
With our choice of $\psi$ and in view of  (\ref{=0}) and (\ref{estmQphi}), 
we see that we can choose $s$ and $R$
$(s < R)$ such that
\begin{equation}\label{2ndestm}
Q(\phi,\phi) < \delta.
\end{equation}
Combining (\ref{1stestm}) and (\ref{2ndestm}) we get 
$Q(\phi_{\epsilon}, \phi_{\epsilon}) < 0$, and the proof of the 
theorem is complete.

\qed

We make a final remark on the following question brought up by the referee: other than tubes over totally geodesic submanifolds, are there  quantum tubes that do not have pure point spectrum below the continuum?

The question is very difficult, even in the surface case. The following is a conjecture people are working on:

{\bf Conjecture.} \begin{it}
Let $\Sigma$ be a complete surface embedded into $\R^3$ such that the second fundamental form goes to zero at infinity. Suppose the Gaussian curvature is integrable. If there  is  no pure point spectrum below $\sigma_{ess}(\Delta)$ then $\Sigma$ must be totally geodesic.\end{it}

Under the additional assumption that the Gaussian curvature is positive, the conjecture is true. This essentially follows from Theorem 1.3 of ~\cite{ll-1}.

\bibliographystyle{abbrv} 
\bibliography{new051007,unp051007,local}

\end{document}